\documentclass[12pt]{article}

\usepackage[table]{xcolor}
\usepackage{amsmath, amsthm}
\usepackage{amsfonts,amssymb}
\usepackage{hyperref}
\usepackage{authblk}
\usepackage[all]{xy}
\usepackage{enumitem}
\usepackage{stmaryrd}
\usepackage[makeroom]{cancel}
\usepackage{bbm}
\usepackage{fancyvrb}

\usepackage[linesnumbered,boxed,vlined]{algorithm2e}

\usepackage{xstring} 

\usepackage{tikz}
\usetikzlibrary{intersections}

\newif\ifextmath
\extmathfalse

\newcommand{\Z}{\mathbb{Z}}
\newcommand{\Q}{\mathbb{Q}}

\newcommand{\C}{\mathbb{C}}
\renewcommand{\L}{\mathcal{L}}
\newcommand{\B}{\mathcal{B}}
\newcommand{\F}{\mathbb{F}}
\renewcommand{\t}{\theta}
\renewcommand{\l}{\lambda}

\newcommand{\sep}{\text{sep}}
\newcommand{\ins}{\text{ins}}
\newcommand{\Tr}{\operatorname{Tr}}

\newcommand{\Hom}{\operatorname{Hom}}
\newcommand{\Aut}{\operatorname{Aut}}
\newcommand{\Id}{\operatorname{Id}}
\newcommand{\caract}{\operatorname{char}}

\newcommand{\floor}[1]{\left\lfloor #1 \right\rfloor}

\SetKw{True}{True}
\SetKw{False}{False}
\newcommand{\la}{\leftarrow}
\SetKw{FAIL}{FAIL}
\SetKw{Goto}{go to line}

\numberwithin{equation}{section}

\newtheorem{thm}[equation]{Theorem}
\newtheorem{lem}[equation]{Lemma}
\newtheorem{pro}[equation]{Proposition}
\newtheorem{cor}[equation]{Corollary}
\theoremstyle{definition}
\newtheorem{de}[equation]{Definition}

\newtheorem{ex}[equation]{Example}

\makeatletter
\newcommand{\subjclass}[2][2020]{%
  \let\@oldtitle\@title%
  \gdef\@title{\@oldtitle\footnotetext{#1 \emph{Mathematics subject classification:} #2}}%
}
\newcommand{\keywords}[1]{%
  \let\@@oldtitle\@title%
  \gdef\@title{\@@oldtitle\footnotetext{\emph{Key words and phrases.} #1.}}%
}
\let\c@table\c@equation
\let\c@figure\c@equation
\makeatother

\makeatletter
\let\oldnl\nl
\newcommand{\nonl}{\renewcommand{\nl}{\let\nl\oldnl}}
\makeatother

\title{Computing the trace of an algebraic point \\ on an elliptic curve}
\author{Nicolas Mascot\thanks{Trinity College Dublin, \href{mailto:mascotn@tcd.ie}{mascotn@tcd.ie}}, Denis Simon\thanks{Universit\'e de Caen-Normandie, \href{mailto:denis.simon@unicaen.fr}{denis.simon@unicaen.fr}}}

\VerbatimFootnotes

\begin{document}

\maketitle

\begin{abstract}
We present a simple and efficient algorithm to compute the sum of the algebraic conjugates of a point on an elliptic curve.
\end{abstract}


\section{Introduction}

Let~$E$ be the elliptic curve over~$\Q$ defined by the equation~$y^2=x^3+x+15$, and let~$\theta$ be a root of the irreducible polynomial~$T(t) = t^3-135t-408 \in \Q[t]$, so that the point~$P = (\frac{\theta}8-1,\frac{\theta^2}{32} - \frac{11\theta}{32} - \frac{19}4)$ lies on~$E$. The sum in~$E$ of the three algebraic conjugates of~$P$ must then lie in~$E(\Q)$, but which point is it exactly?

A computer package such as \cite{gp} reveals that the images of~$P$ under the three embeddings of~$\Q(\theta)$ into~$\C$ are
\[ (-2.202\cdots,1.451\cdots), (-1.410\cdots,-3.283\cdots), (0.613\cdots,-3.980\cdots), \]
which add up in~$E$ to~$(2.000\cdots,-5.000\cdots)$, so the answer must be~$(2,-5)$. However, this approach is not rigorous since it involves approximations. Besides, if we had been working over another ground field than~$\Q$, such as~$\Q(t)$ for example, it would not have been so easy to find an analogue of~$\C$ in which explicit computations are easy.

Another approach would consist in working in a splitting field of~$T(t)$; however, this approach is ineffective since in general, the degree of the splitting field of a polynomial of degree~$n$ may be as large as~$n!$, which is prohibitive even for moderate values of~$n$.

The purpose of this article is to present a simple and efficient algorithm to compute the sum of the algebraic conjugates of a point on an elliptic curve and which is valid over any ground field.

Our implementation of these algorithms has been available since March 2022 in the development version of \cite{gp} under the name \verb?elltrace?.

\section{Formal definition and first properties}

Fix a field~$K$ and an elliptic curve~$E$ over~$K$.

\begin{de}\label{de:trace}
Let~$\Omega$ be an algebraic closure of~$K$, and let~$i \in \Hom(K,\Omega)$ be an embedding. Let~$L$ be a finite extension of~$K$ of inseparable degree~$[L:K]_\ins$, and let~$P \in E(L)$. The \emph{trace}~$\Tr(L/K,i,P)$ of~$P$ relative to the extension~$L/K$ and to the embedding~$i$ is the point~$T \in E(K)$ such that
\begin{equation} T^i = [L:K]_\ins \sum_{j \in \Hom_i(L,\Omega)} P^j, \label{eqn:deftrace} \tag{$\dagger$} \end{equation}
where the sum is meant in the sense of the group law of~$E$ and ranges over the embeddings~$j$ of~$L$ into~$\Omega$ which extend~$i$. We will prove in Proposition~\ref{prop:welldef} below that the right hand side of \eqref{eqn:deftrace} indeed lies in $E\big(i(K)\big)$.
\end{de}

\begin{ex}
Suppose that~$L = K(\theta)$ for some~$\theta \in L$. Let
\[ \prod_{n=1}^{[L:K]_\sep} (t-\theta_n)^{[L:K]_\ins} \in \Omega[t] \]
be the factorisation over~$\Omega$ of the image by $i$ of the minimal polynomial of~$\theta$ over~$K$, and write~$P = \big(x(\theta),y(\theta)\big)$ where~$x(t),y(t) \in K[t]$. Then
\[ \Tr(L/K,i,P) = i^{-1}\left( [L:K]_\ins\sum_{n=1}^{[L:K]_\sep} \big(x(\theta_n),y(\theta_n)\big) \right). \]
\end{ex}

\begin{lem}\label{lem:sep_pow_p}
Let~$K$ be a field of characteristic~$p > 0$, and let~$L$ be a finite extension of~$K$ of inseparable degree~$[L:K]_\ins = p^d$. Then for all~$\alpha \in L$,~$\alpha^{p^d}$ is separable over~$K$.
\end{lem}

\begin{proof}
Let~$A(t) \in K[t]$ be the minimal polynomial of~$\alpha$, and let~$e \geq 0$ be the largest integer such that~$A(t) \in K[t^{p^e}]$. Then~$A(t) = B(t^{p^{e}})$ for some separable~$B(t) \in K[t]$, so~$\alpha^{p^e}$ is separable over~$K$. Besides,
\[ p^e = [K(\alpha):K]_\ins \mid [L:K]_\ins = p^d, \]
so~$e \leq d$.
\end{proof}

\begin{pro}\label{prop:welldef}
Definition~\ref{de:trace} makes sense.
\end{pro}

\begin{proof}
Let~$Q = [L:K]_\ins P$. We must prove that the point
\[ R = \sum_{j \in \Hom_i(L,\Omega)} Q^j, \]
which \emph{a priori} merely lies in~$E(\Omega)$, actually lies in~$E\big(i(K)\big)$.

Identifying~$K$ with~$i(K)$, we may assume without loss of generality that~$K \subseteq \Omega$ and that~$i = \Id$.

We first prove that  for all~$j \in \Hom_i(L,\Omega)$, the point~$Q^j$ is defined over the separable closure~$K^\sep$ of~$K$ in~$\Omega$. If~$\caract K = 0$ there is nothing to prove, so assume that~$\caract K = p > 0$. Let~$d \in \Z_{\geq 0}$ be such that~$[L:K]_\ins = p^d$. Recall \cite[p.~320]{DS} that the multiplication-by-$p$ map~$[p]$ on~$E$ factors through the Frobenius map~$(x,y) \mapsto (x^p,y^p)$; Lemma~\ref{lem:sep_pow_p} thus implies that~$Q^j = [L:K]_\ins P^j = [p]^d P^j$ lies in~$E(K^\sep)$, as desired.

As a result,~$R$ is an element of~$E(K^\sep)$, which is fixed by~$\Aut_K(K^\sep)$ since~$\Aut_K(K^\sep)$ permutes the elements of~$\Hom_K(L,K^\sep)$, and hence actually an element of~$E(K)$.
\end{proof}

\begin{lem}\label{lem:tr_indep_i}
The trace~$\Tr(L/K,i,P)$ does not depend on~$\Omega$ nor on~$i$.
\end{lem}

\begin{proof}
Let~$i'$ be another embedding of~$K$ into another algebraic closure~$\Omega'$ of~$K$. There exists at least one morphism~$\varphi : \Omega \longrightarrow \Omega'$ such that~$i' = \varphi i$; it is then clear that the maps
\[ \begin{array}{ccc} \Hom_i(L,\Omega) & \longleftrightarrow & \Hom_{i'}(L,\Omega') \\ j & \longmapsto & \varphi j \\ \varphi^{-1} j' & \longmapsfrom & j' \end{array} \]
are bijections that are inverses of each other. Therefore
\[ \Tr(L/K,i',P)^{i'} = [L:K]_\ins \hspace{-5mm} \sum_{j' \in \Hom_{i'}(L,\Omega')} \hspace{-5mm} P^{j'} = \left([L:K]_\ins \hspace{-5mm} \sum_{j \in \Hom_{i}(L,\Omega)} \hspace{-5mm} P^{j}\right)^\varphi = \big(\Tr(L/K,i,P)^i\big)^{\varphi}, \] 
as wanted.
\end{proof}

Lemma~\ref{lem:tr_indep_i} justifies the notation~$\Tr^L_K P$ without reference to $i$, which we adopt from now on.

\begin{pro} \label{prop:morph}
$\Tr^L_K$ is a group morphism from~$E(L)$ to~$E(K)$.
\end{pro}

\begin{proof}
Obvious from Definition~\ref{de:trace}.
\end{proof}

\begin{thm}\label{thm:trans}
Suppose we have a tower~$K \subseteq L \subseteq M$ of finite extensions and a point~$P \in E(M)$. Then
\[ \Tr^M_K P = \Tr^L_K \Tr^M_L P. \]
\end{thm}

\begin{proof}
We compute that
\begin{align*}
\Tr(M/K,i,P)^i &= [M:K]_\ins \sum_{k \in \Hom_i(M,\Omega)} P^k \\
&= [L:K]_\ins \sum_{j \in \Hom_i(L,\Omega)} [M:L]_\ins \sum_{k \in \Hom_j(M,\Omega)} P^k \\
&= [L:K]_\ins \sum_{j \in \Hom_i(L,\Omega)} \Tr(M/L,j,P)^j \\
& = [L:K]_\ins \sum_{j \in \Hom_i(L,\Omega)} (\Tr^M_L P)^j \\
&= \Tr(L/K,i,\Tr^M_L P)^i. \qedhere \\
\end{align*}
\end{proof}

\begin{cor} \label{cor:subdef}
If~$K \subseteq L \subseteq M$ and if~$P \in E(L)$, then 
\[ \Tr^M_K P = [M:L] \Tr^L_K P. \]
\end{cor}

\begin{proof}
It s clear from Definition~\ref{de:trace} that 
\[ \Tr^M_L P = [M:L]_\ins [M:L]_\sep P = [M:L] P,\]
so the result follows from Theorem~\ref{thm:trans} and Proposition~\ref{prop:morph}.
\end{proof}

\pagebreak

\section{The algorithm}

We are now going to describe an algorithm to compute the trace of a point, which we may assume is not the point at infinity~$O$ in view of Proposition~\ref{prop:morph}.

\subsection{Assumptions and prerequisites}

We suppose that the ground field~$K$ is such that its elements may be represented exactly on a computer, and that algorithms are known for the four field operations in~$K$. For example,~$K$ could be a finite field, a number field, or a field of rational fractions in finitely many variables over such a field.

It is then possible to perform computations in extensions of~$K$ of the form~$K[t]/T(t)$ where~$T(t) \in K[t]$ is irreducible, and thus in simple extensions of~$K$.

It is also possible to perform linear algebra computations over~$K$, and, in particular, given a matrix~$M$ with coefficients in~$K$, to compute a \emph{kernel matrix}~$Z$ of~$M$ in \emph{echelon form}, that is to say a matrix whose columns represent a basis of the (right) kernel of~$M$, and which is such that the number of consecutive zeroes at the bottom of the~$j$-th column of~$Z$ is a strictly decreasing function of~$j$. This is the convention followed by the \cite{gp} function \verb?matker?.

In particular, this makes it possible to determine the minimal polynomial over $K$ of an algebraic element.

We restrict ourselves to the case of a point defined over a simple extension of~$K$, that is to say an extension $K(\theta)/K$ generated by a single element $\theta$, as we expect this to cover most practical applications. If needed, the general case can be handled by a repeated application of Theorem~\ref{thm:trans}.

\subsection{The separable case}

We begin by the simpler case where the extension is separable.  We will show in the next section how to reduce to this case.

\pagebreak

\begin{algorithm}[H]
\nonl \TitleOfAlgo{EllTraceSep}
\KwIn{An elliptic curve~$E:y^2+a_1xy+a_3y=x^3+a_2x^2+a_4x+a_6$ with the~$a_i \in K$, an irreducible \emph{and separable} polynomial~$T(t) \in K[t]$, and a pair  of polynomials~$x_P(t),y_P(t) \in K[t]$ of degree~$< \deg T$ such that~$P=\big(x_P(\theta),y_P(\theta)\big) \in E\big(K(\theta)\big)$ where~$\theta$ is the class of~$t$ in~$K[t]/T(t)$.}
\KwOut{The trace~$\Tr^{K(\theta)}_K P \in E(K)$.}
$d \la \deg T$\;
\If{$x_P(t)$ is constant}
{
\eIf{$y_P(t)$ is also constant}
{\Return~$[d]P$\;}
{\Return $O$\;} \label{alg:elltrace_sep:trivial}
}
$\L \la$ a vector of length~$d+1$\;
$\L[1] \la 1$;~$\L[2] \la x_P(\theta)$;~$\L[3] \la y_P(\theta)$\;
\For{$j \la 4$ \KwTo $d+1$}
{$\L[j] = x_P(\theta) \L[j-2]$\;}
$M \la$ a matrix of size~$d \times (d+1)$\;
\For{$j \la 1$ \KwTo $d+1$}
{$j$-th column of~$M \la$ coefficients of~$\L[j]$ wrt.\ $1,\theta,\theta^2,\cdots,\theta^{d-1}$\;}
$Z \la$ leftmost column of a kernel matrix of~$M$ in echelon form\;
$U(x) \la Z[1] + \sum_{j=1}^{\floor{\frac{d+1}2}} Z[2j] x^j$;~$V(x) \la \sum_{j=0}^{\floor{\frac{d-2}2}} Z[2j+3] x^j$\; \label{alg:elltrace_sep:fP}
\If{$V(x)$ is the $0$ polynomial}{\Return $O$\;} \label{alg:elltrace_sep:V0}
$X(x) \la$ minimal polynomial of $x_P(\theta)$ over $K$\;
$R(x) \la (x^3+a_2x^2+a_4x+a_6) V(x)^2 + (a_1 x + a_3) U(x) V(x) - U(x)^2$\;
$S(x) \la R(x) / X(x)$\;
\If{$\deg S(x) = 0$}
{\Return $O$\; \label{alg:elltrace_sep:Rconst}} 
$x_Q \la$ unique root of~$S(x)$\; \label{alg:elltrace_sep:xQ}
$y_Q \la -U(x_Q)/V(x_Q)$\; \label{alg:elltrace_sep:yQ} 
$Q \la (x_Q,y_Q)$\;
$d_P \la$ largest~$j \leqslant d$ such that~$Z[j+1] \neq 0$\;  \label{alg:elltrace_sep:dP}
\Return~$[-d/d_P] Q$\;
\caption{Trace of an algebraic point, separable case}
\label{alg:elltrace_sep}
\end{algorithm}

\pagebreak

\begin{proof}
Let~$L = K(\theta)$, and let~$\Omega \supseteq K$ be an algebraic closure of~$K$. In view of Proposition~\ref{prop:welldef}, we may compute~$\Tr^L_K P$ as~$\Tr(L/K,\Id,P)$.

Let~$K(P) = K\big(x_P(\theta),y_P(\theta)\big) \subseteq L$ be the field of definition of~$P$, and let~$d_P = [K(P) : K]$, so that~$d_P \mid d = [L:K]$. We will prove below that line~\ref{alg:elltrace_sep:dP} correctly determines~$d_P$. Let~$Q' = \Tr^{K(P)}_K P$, and let~$Q = [-1] Q'$, so that 
\begin{equation} \Tr^L_K P = [L:K(P)] Q' = [-d/d_P] Q \tag{$\star$} \label{eqn:TrKPdP} \end{equation}
by Corollary~\ref{cor:subdef}.

By assumption,~$L$, and therefore~$K(P)$, are separable over~$K$; therefore the~$d_P$ points~$P^j$,~$j \in \Hom_K(K(P),\Omega)$, are pairwise distinct.

Suppose that $x_P(\theta) \in K$. If we also have $y_P(\theta) \in K$, then $P \in E(K)$, so trivially $\Tr^L_K P = [d] P$ by Definition \ref{de:trace}, whereas if $y_P(\theta) \not \in K$, then necessarily $d_P=2$ and the two $P^j$, which are distinct, are thus negatives of each other, so $\Tr^{K(P)}_K P = O$. This justifies the lines up to line~\ref{alg:elltrace_sep:trivial}.

We assume that~$x_P(\theta) \not \in K$ from now on. The divisor 
\[ D = \sum_{j \in \Hom_K(K(P),\Omega)} P^j + Q - (d_P + 1) O \]
is defined over~$K$ by separability, and is principal by definition of~$Q$. Let~$f_P \in K(E)$ have divisor~$D$; its degree is~$d_P+1$ if~$Q \neq O$, and is~$d_P$ if~$Q = O$. Conversely, any function~$f \in K(E)$ which vanishes at all the~$P^j$ satisfies~$\deg f \geq \deg f_P$. Thus~$f_P$ is a function of \emph{minimal degree} among those that vanish at all of the~$P^j$. In other words,~$f_P$ has minimal degree among the functions that are defined over~$K$ and vanish at~$P$, and any element of $K(E)$ which has the same degree as $f_P$, vanishes at $P$, and has no other pole than $O$ must be a scalar multiple of $f_P$.

Since~$d_P \leq d$, the function~$f_P$ lies in the Riemann-Roch space attached to the divisor~$(d+1)O$, which has dimension~$d+1$ and admits the~$K$-basis
\[ \B = \left\{ \begin{array}{ll} \{1,x,y,x^2,xy,\cdots,x^{\frac{d}2},x^{\frac{d}2-1}y \}, & d \text{ even}, \\ \{1,x,y,x^2,xy,\cdots,x^{\frac{d-1}2-1}y,x^{\frac{d+1}2} \}, & d \text{ odd}. \end{array} \right. \]

Number the elements of~$\B$ as~$b_1=1$,~$b_2=x$,~$b_3=y$,~$b_4=x^2$, and so on, so that~$\deg b_j = j$ for all~$j \geq 2$. Besides,~$b_{j} = x b_{j-2}$ for all~$4 \leq j \leq d+1$, so~$\L$ is the vector of values in~$L$ of the elements of~$\B$ at~$P$. The algorithm then expresses these values on the~$K$-basis~$1,\theta,\cdots,\theta^{d-1}$ of~$L$, whence the matrix~$M$ whose kernel consists of the coordinates on~$\B$ of the functions in the~$K$-span of~$\B$ that vanish at~$P$. As~$\deg b_j$ is a strictly increasing function of~$j$, the leftmost column in a kernel matrix of~$M$ in echelon form contains the coordinates on~$\B$ of some nonzero scalar multiple of~$f_P$, which we may assume is exactly~$f_P$ since~$f_P$ is only defined up to scaling. We express~$f_P$ as~$U(x)+V(x)y$ at line~\ref{alg:elltrace_sep:fP}.

$\bullet$ Suppose for now that~$V(x)=0$. We claim that $2P \neq O$ in this case. Indeed, if we had~$2P=O$,  then the~$2$-division polynomial
\[ \psi_2 = 2y+a_1x+a_3 \in K(E)^\times \] would vanish at $P$. If $\caract K = 2$, this would contradict our assumption that~$x_P(\theta) \not \in K$. If $\caract K \neq 2$, then $\psi_2$ cannot be a scalar multiple of $f_P$ since $V(x)=0$, so we would have~$\deg f_P < \deg \psi_2 = 3$ whence~$\deg f_P=2$, so~$f_P \in K(E)$ would be a nonzero scalar multiple of $x-x_P(\theta)$, which would again contradict our assumption that $x_P(\theta) \not \in K$.

As~$f_p = U(x)$ is invariant by the elliptic involution~$[-1]$, so is its divisor of zeroes. As~$P \not \in E(K)$ whereas~$Q \in E(K)$, it follows that for each~$j \in \Hom_K(K(P),\Omega)$, there exists~$j' \in \Hom_K(K(P),\Omega)$ such that~$P^{j'} = [-1] P^j$. As we cannot have~$P^{j'} = P^j$ since~$2P \neq O$, the $P^j$ can be partitioned into pairs of opposites, so $\Tr^{K(P)}_K P = O$ by Definition \ref{de:trace}.  Theorem \ref{thm:trans} thus shows that $\Tr^L_K P = O$, so we return $O$ at line~\ref{alg:elltrace_sep:V0}.

$\bullet$ Suppose from now on that $V(x) \neq 0$. In this case, we claim that $x_P(\theta)$ must generate~$K(P)$ over $K$. For else its minimal polynomial $X(x)$ over $K$ would be of degree $\leq d_P/2$ in $x$, and thus of degree $\leq d_P$ as a function on~$E$, so that $f_P$ would be a nonzero multiple of $X(x)$ by minimality of $\deg f_P$, contradicting our assumption that $V(x)\neq 0$. In particular, it follows that 
\[ X(x) = \prod_{j \in \Hom(K(P),\Omega)} \big(x-x_P(\theta^j)\big). \]

Let~$f'_P = f_P \circ [-1] = U(x) + V(x) (-y-a_1x-a_3)  \in K(E)$. Then 
\begin{align*}
-f_P f'_P &= \big(V(x)y+U(x)\big)\big(V(x)y-U(x)\big) + \big(V(x)y+U(x)\big)V(x)(a_1x+a_3) \\
&= V(x)^2y^2-U(x)^2 + V(x)^2 (a_1x+a_3)y + U(x) V(x) (a_1x+a_3) \\
&= V(x)^2(x^3+a_2x^2+a_4x+a_6) - U(x)^2 + U(x) V(x) (a_1x+a_3) \\
& = R(x).
\end{align*}
Therefore the divisor of~$R(x)$ is
\begin{align*}
 \big(R(x)\big) &= D + [-1]^\times D \\
& = \sum_{j \in \Hom(K(P),\Omega)} (P^j + [-1]P^j) + Q + Q' -2(d_P+1)O \\
&= \left( \prod_{j \in \Hom(K(P),\Omega)} \big(x-x_P(\theta^j)\big) \right) + Q + Q' - 2O \\
&= \big(X(x)\big) + Q + Q' - 2O.
\end{align*}
It follows that~$S(x) = R(x)/X(x)$ is constant if and only if $Q=O$. In particular, in view of \eqref{eqn:TrKPdP}, it is legitimate to return~$O$ at line~\ref{alg:elltrace_sep:Rconst}.

Suppose from now on that~$Q\neq O$, and write $Q=(x_Q,y_Q)$. Then $S(x)$ is a nonzero scalar multiple of $x-x_Q$ , so we correctly recover~$x_Q$ at line~\ref{alg:elltrace_sep:xQ}. If we had~$V(x_Q)=0$, then as~$f_P = U(x)+V(x)y$ vanishes at~$Q$, we would also have~$U(x_Q) = 0$, so we could write~$f_P = (x-x_Q) g_P$ for some~$g_P \in K(E)$. But as $x_P(\theta) \not \in K$, $g_P$ would vanish at all the~$P^j$, which would violate the minimality of~$\deg f_P$. Therefore~$V(x_Q)\neq 0$, so we are able to recover~$y_Q$ at line~\ref{alg:elltrace_sep:yQ}.

Finally, we determine~$d_P$ at line~\ref{alg:elltrace_sep:dP} from the fact that~$\deg b_j = j$ for all~$j \neq 1$ and that~$\deg f_P = d_P+1$ since~$Q \neq O$, which allows us to determine~$\Tr^L_K P$ by \eqref{eqn:TrKPdP}.
\end{proof}

The complexity of this algorithm is dominated by the cost of writing down the list $\L$, of determining a kernel matrix of $M$ in echelon form, and of determining the minimal polynomial $X(x)$. Each of these steps requires~$O(d^3)$ operations in $K$ if one uses naive algorithms for arithmetic in $K[t]$ and linear algebra over $K$. This complexity can be reduced to $O\big(dM(d)+d^\omega\big)$ if one uses fast algorithms for polynomial arithmetic and linear algebra, where polynomials of degree $d$ can be multiplied in $O\big(M(d)\big)$ and $\omega$ is a complexity exponent for linear algebra in dimension $d$.

\pagebreak

\begin{ex}
Let us demonstrate the execution of the algorithm on the example in the introduction. In this example, we have $K=\Q$, $L = \Q(\theta)$, $E$ is given by $y^2=x^3+x+15$, and $P=\left(\frac18 \theta-1,\frac1{32}\theta^2-\frac{11}{32}\theta-\frac{19}4\right)$ where $\theta$ has minimal polynomial $T(t) = t^3-135t-408 \in \Q[t]$.

We have $d = \deg T = 3$. Since $x_P(\theta) \not \in \Q$, we begin by computing
\[ \L = [ 1, x_P(\theta), y_P(\theta), x_P(\theta)^2 ] = \left[ 1, \frac18 \theta-1,\frac1{32}\theta^2-\frac{11}{32}\theta-\frac{19}4, \frac1{64}\theta^2-\frac{1}{4}\theta+1 \right]. \]
Thus $M$ is the matrix
\[ M = \left[ \begin{matrix} 1 & -1 & -\frac{19}4 & 1 \\ 0 & \frac18 & -\frac{11}{32} & -\frac14 \\ 0 & 0 & \frac1{32} & \frac1{64} \end{matrix} \right]. \]

We find that $\ker M$ has dimension $1$ and is spanned by $[-22,5,-4,8]$, so we may take $f_P = -22+5x-4y+8x^2 = U(x)+V(x)y$ where $U(x)=8x^2+5x-22$ and $V(x)=-4$.

As $V(x) \neq 0$, we determine that the minimal polynomial of $x_P(\theta)$ over $\Q$ is
\[ X(x) = x^3 + 3x^2 + \frac{57}{64}x - \frac{61}{32}, \]
and we compute that
\[ R(x) = (x^3+x+1)V(x)^2 - U(x)^2 = -64x^4 - 64x^3 + 327x^2 + 236x - 244. \]
This leads to $S(x) = R(x)/X(x) = -64x+128 = -64(x-2)$, so $x_Q=2$, and $y_Q = -U(x_Q)/V(x_Q) = 5$, meaning that $Q = (2,5)$.

Finally, since $d_P +1 = \deg f_P = 4$, we have $d_P = 3$, whence as expected
\[ \Tr^{\Q(\theta)}_\Q P = [-1]Q = (2,-5) \in E(\Q). \]
\end{ex}

\begin{ex}
As a second example, let $E$ be the elliptic curve $y^2=x^3+x^2+1$ over $K=\F_3$, let $L = \F_{3^6} = K(\theta)$ where $\theta$ is a root of $T(t) = t^6+t^5+t^4+t^3+t^2+t+1 \in K[t]$, and let
\[ P = ( \t^5+\t^2, \t^4+\t^3+2) \in E(L). \]
This time we have $d=6$, so we compute
\[ \L =  [ 1, x_P(\theta), y_P(\theta), x_P(\theta)^2, x_P(\theta)y_P(\theta), x_P(\theta)^3, x_P(\theta)^2 y_P(\theta) ] \] 
which leads to
\[ M = \left[\begin{matrix}
1 & 0 & 2 & 2 & 2 & 2 & 2\\
0 & 0 & 0 & 0 & 0 & 0 & 0\\
0 & 1 & 0 & 0 & 2 & 2 & 2\\
0 & 0 & 1 & 1 & 2 & 2 & 0\\
0 & 0 & 1 & 1 & 2 & 2 & 0\\
0 & 1 & 0 & 0 & 2 & 2 & 2\\
\end{matrix} \right] .\]
We find that $\dim \ker M = 4$, and that
\[ \left[\begin{matrix}
0 & 2 & 2 & 1\\
0 & 1 & 1 & 1\\
2 & 1 & 1 & 0\\
1 & 0 & 0 & 0\\
0 & 1 & 0 & 0\\
0 & 0 & 1 & 0\\
0 & 0 & 0 & 1\\
\end{matrix} \right] \]
is a kernel matrix of $M$ in echelon form, so we take $f_P = 2y+x^2$, which corresponds to $U(x)=x^2$ and $V(x)=2$.

As $V(x) \neq 0$, we determine that the minimal polynomial of $x_P(\t)$ over~$K$ is
\[ X(x) = x^3 + x^2 + x + 2, \]
and we compute that
\[ R(x) = (x^3+x^2+1)V(x)^2-U(x)^2 = 2x^4 + x^3 + x^2 + 1 = 2 (x+1) X(x). \]
Therefore $x_Q = 2$, and $y_Q = -U(x_Q)/V(x_Q) = 1$.

Besides, we have $d_P = \deg f_P -1 = 3$, so this time $K(P) = \F_{3^3} \subsetneq L$. We conclude that
\[ \Tr^{\F_{3^6}}_{\F_3} P = [-2] Q = (2,1). \]
\end{ex}

\begin{ex}
As a last example, take  $K=\F_2(\l)$ where $\l$ is an indeterminate, consider the elliptic curve $E$ of equation \[ y^2+\l xy+y=x^3+\l x^2 + \l x, \]
and let
\[ P = (\t^4+\l \t^2-\t+\l, \l \t^4+\t^3+\l^2 \t^2 + \l^2+1) \in E(L), \]
where $L = K(\t)$ and $\t$ is a root of the irreducible polynomial 
\[ T(t) = t^5+\l t^3 + \l t + \l \in K[t]. \]

We find that
\[ M = \left[ \begin{matrix}
1 &  \l &  \l^2 + 1  &    \l^2 &  \l^3 +  \l^2 +  \l  &    \l^2 \\
0 & 1   &    0  &    \l^2   &     \l^3 + 1&  \l^3 +  \l \\
0 &  \l   &   \l^2 &  \l^2 + 1  &          \l^3  &    \l^3 \\
0 & 0  &     1    &    \l      &      \l^2 &  \l^2 + 1 \\
0 & 1  &      \l    &    \l      &      \l^2 &  \l^2 +  \l \\
\end{matrix} \right] \]
has kernel of dimension $1$ which is spanned by
\[ [  \l^3+ \l^2+1,  \l^3+ \l^2, 1,  \l^4+ \l^2, \l^3+ \l^2+ \l, \l^2+1 ], \]
whence $U(x) = ( \l^2 + 1)x^3 + ( \l^4 +  \l^2)x^2 + ( \l^3 +  \l^2)x + ( \l^3 +  \l^2 + 1)$ and $V(x)=( \l^3 +  \l^2 +  \l)x + 1$.

This leads to
\[ R(x) = (x^3+ \l x^2+ \l x)V(x)^2+( \l x+1)U(x) V(x) + U(x)^2 = X(x) \big( ( \l^4 + 1)x + ( \l^4 +  \l^3 +  \l) \big), \]
so $x_Q = \frac{ \l^4 +  \l^3 +  \l}{ \l^4+1}$ and $y_Q = U(x_Q)/V(x_Q) = \frac{ \l^7 + \l^5 + \l^4 + \l^2 + 1}{ \l^6 +  \l^4+ \l^2 + 1}$.

As $d_P = \deg f_P -1 = 5$, we conclude that
\[ \Tr^{\F_2(\l,\t)}_{\F_2(\l)} P = [-1]Q = \left( \frac{ \l^4 +  \l^3 +  \l}{ \l^4+1}, \frac{ \l^2}{ \l^6 +  \l^4+ \l^2 + 1} \right). \]

Note that $T(t)$ happens to have Galois group $S_5$ over $K$, so an approach based on the explicit construction of its splitting field would have been extremely laborious.
\end{ex}

\pagebreak

\subsection{The general case}

We now present an algorithm to handle the general case, when~$T(t)$ is not necessarily separable.

\bigskip

\begin{algorithm}[H]
\nonl \TitleOfAlgo{EllTrace}
\KwIn{An elliptic curve~$E:y^2+a_1xy+a_3y=x^3+a_2x^2+a_4x+a_6$ with the~$a_i \in K$, an irreducible polynomial~$T(t) \in K[t]$, and a pair of polynomials~$x_P(t),y_P(t) \in K[t]$ of degree~$< \deg T$ such that~$P=\big(x_P(\theta),y_P(\theta)\big) \in E\big(K(\theta)\big)$ where~$\theta$ is the class of~$t$ in~$K[t]/T(t)$.}
\KwOut{The trace~$\Tr^{K(\theta)}_K P \in E(K)$.}
$p \la \caract K$\;
\eIf{$p \neq 0$}{
$d \la$ largest integer~$\geq 0$ such that~$T(t) \in K[t^{p^d}]$\;
$S(t) \la$ polynomial such that~$T(t) = S(t^{p^d})$\;
$Q \la [p^d] P$\; \label{alg:elltrace:mul_pd}
\If{$Q=O$}
{\Return~$O$\;}
$x_Q(t), y_Q(t) \la$  polynomials of degree~$< \deg S$ such that~$Q=\big(x_Q(\theta^{p^d}),y_Q(\theta^{p^d})\big)$\;\label{alg:elltrace_xQ}
\Return EllTraceSep$\big(E,S,(x_Q,y_Q)\big)$\;
}
{
\Return EllTraceSep$\big(E,T,(x_P,y_P)\big)$\;
}
\caption{Trace of an algebraic point}
\label{alg:elltrace}
\end{algorithm}

\begin{proof}
(Compare with the proof of Proposition~\ref{prop:welldef}) Let~$L=K[t]/T(t)$, and let~$K \subseteq L_\sep \subseteq L$ be the separable closure of~$K$ in~$L$.

As~$L/L_\sep$ is purely inseparable, the point $Q = [p^d] P$ satisfies
\[ Q = [L:K]_\ins P = [L:L_\sep] P = \Tr^L_{L_\sep} P \in E({L_\sep}) \]
by Definition~\ref{de:trace}. Theorem~\ref{thm:trans} thus ensures that 
\[ \Tr^L_K P = \Tr^{L_\sep}_K \Tr^L_{L_\sep} P = \Tr^{L_\sep}_K Q. \]
Finally, observe that while the~$\theta^n$ for~$0 \leq n < \deg T$ form a~$K$-basis of~$L$, the~$(\theta^{p^d})^n$ for~$0 \leq n < \deg S$ form a~$K$-basis of~$L_\sep$; therefore the polynomials~$x_Q(t)$ and~$y_Q(t)$ introduced at line~\ref{alg:elltrace_xQ} exist and are unique, and can be read directly off the coordinates of~$Q$. 
\end{proof}

The complexity of this algorithm is dominated by the multiplication by~$p^d$ at line \ref{alg:elltrace:mul_pd} and by the execution of EllTraceSep, and is therefore bounded by~$O([L:K]^2 \log [L:K]_\ins + [L:K]_\sep^3)$ operations in $K$, and by 
\[ O\big(M([L:K]) \log [L:K]_\ins + [L:K]_\sep M([L:K]_\sep) + [L:K]_\sep^\omega\big) \]
if one uses fast algorithms for polynomial arithmetic and for linear algebra.

\begin{ex}
Let $E : y^2+xy+\l y = x^3+x^2$ over $K = \F_2(\l)$ where $\l$ is an indeterminate, let $L=K(\theta)$ where $\theta$ is a root of the irreducible polynomial 
\[ T(t) = t^4+t^2+\l^4+\l^3 \in K[t], \]
and let
\[ P = \big( \t^2+\t, \t^3+(\l+1) \t  + \l^2+\l \big) \in E(L). \]

We have $p=2$, $d=1$, $[L:K]=4$, $[L:K]_\sep = [L:K]_\ins = 2$, and~$L_\sep = K(\theta^2)$ where the minimal polynomial of $\theta^2$ is $S(t) =  t^2+t+\l^4+\l^3$.

We calculate
\[ Q = [2] P = \big(\l^4 + \l^3 + \l^2 + \l + 1, (\l^4 + \l^3 + \l^2 + 1) \t^2 + \l^5 + \l \big) \]
which lies in $E(L_\sep)$ as expected, so we call EllTraceSep with 
\[ x_Q(t) = \l^4 + \l^3 + \l^2 + \l + 1, \ y_Q(t) = (\l^4 + \l^3 + \l^2 + 1) t + \l^5 + \l. \]
As $x_Q(t)$ is actually constant as a polynomial in $t$ whereas $y_Q$ is not, \linebreak EllTraceSep immediately reports that $\Tr^{L_\sep}_K Q = O$, so we conclude that
\[ \Tr^{\F_2(\l,\t)}_{\F_2(\l)} P = O. \]
\end{ex}

\renewcommand{\abstractname}{Acknowledgements}
\begin{abstract}
We thank Christophe Delaunay for having asked, in relation to his work on \cite{Christophe}, how to efficiently compute the trace of a point on an elliptic curve; and Bill Allombert for having relayed this question during an informal online \cite{gp} workshop.

Our work on the answer provided in this article was initiated during the January 2022 \cite{gp} workshop, and we express our gratitude to the organisers of this conference.
\end{abstract}


\begin{thebibliography}{PARI/GP}
\addcontentsline{toc}{section}{References}

\bibitem[BBD21]{Christophe} Francesco Battistoni, Sandro Bettin, and Christophe Delaunay, \textbf{Ranks of elliptic curves over~$\Q(T)$ of small degree in~$T$}, arXiv preprint \href{https://arxiv.org/abs/2109.00738}{2109.00738}.

\bibitem[DS05]{DS} Diamond, Fred; Shurman, Jerry,
\textbf{A First Course in Modular Forms}.
Graduate Texts in Mathematics, 228. Springer-Verlag, New York, 2005.

\bibitem[PARI/GP]{gp}The PARI~Group, \textbf{PARI/GP} version \texttt{2.14.0}, Univ. Bordeaux, 2022, \url{http://pari.math.u-bordeaux.fr/}.

\end{thebibliography}
\end{document}